\documentclass[11pt]{article}

\usepackage{amssymb, amsfonts, amsmath, amsthm, mathabx}
\usepackage{lineno}

\newtheorem{theorem}{Theorem}[section]
\newtheorem{lemma}[theorem]{Lemma}
\newtheorem{proposition}[theorem]{Proposition}

\theoremstyle{definition}
\newtheorem{definition}[theorem]{Definition}
\newtheorem{notation}[theorem]{Notation}
\newtheorem{remark}[theorem]{Remark}
\newtheorem{example}[theorem]{Example}
\newtheorem{review}[theorem]{Review}

\newcommand{\cA}{ {\cal A} }

\newcommand{\bC}{ {\mathbb C} }

\newcommand{\fI}{ {\mathcal I} }

\newcommand{\bN}{ {\mathbb N} }
\newcommand{\fO}{ {\mathcal O} }

\newcommand{\cS}{{\cal S}}
\newcommand{\fU}{ {\mathcal U} }
\newcommand{\fV}{ {\mathcal V} }

\newcommand{\cat}{ \mbox{Cat} }
\newcommand{\kk}{ \kappa }

\newcommand{\piciup}{ \widecheck{\pi} }
\newcommand{\picirc}{ \stackrel{\circ}{\pi} }
\newcommand{\rhociup}{ \widecheck{\rho} }
\newcommand{\rhocirc}{ \stackrel{\circ}{\rho} }
\newcommand{\thetaciup}{ \widecheck{\theta} }

\begin{document}

$\ $

\begin{center}
{\bf\Large Non-crossing linked partitions, the partial}

\vspace{6pt}

{\bf\Large order \boldmath{$\ll$} on \boldmath{$NC(n)$}, and the 
           S-transform}

\vspace{14pt}

{\large Alexandru Nica \footnote{Research supported by a 
Discovery Grant from NSERC, Canada.} }

\vspace{10pt}

\end{center}

\begin{abstract}

\noindent
The paper establishes a connection between two recent 
combinatorial developments in free probability: the non-crossing 
linked partitions introduced by Dykema in 2007 to study the 
S-transform, and the partial order $\ll$ on $NC(n)$ introduced
by Belinschi and Nica in 2008 in order to study relations
between free and Boolean probability. More precisely, one has 
a canonical bijection between $NCL(n)$ (the set of all non-crossing 
linked partitions of $\{ 1, \ldots , n \}$) and the set
$\{ ( \alpha , \beta ) \mid \alpha , \beta \in NC(n), \
\alpha \ll \beta \}$. As a consequence of this bijection, one 
gets an alternative description of Dykema's formula expressing the
moments of a noncommutative random variable $a$ in terms of the 
coefficients of the reciprocal $S$-transform $1/S_a$. Moreover, due
to the Boolean features of $\ll$, this formula can be simplified 
to a form which resembles the moment-cumulant formula from c-free 
probability.
\end{abstract}

\vspace{6pt}

\begin{center}
{\bf\large 1. Introduction}
\end{center}
\setcounter{section}{1}

\noindent
This note puts into evidence a connection between two seemingly 
unrelated combinatorial developments in free probability, which 
appeared recently in the papers \cite{BN08} and \cite{D07}.

The basic structure for combinatorial considerations in free 
probability is the poset $( NC(n), \leq)$, where $NC(n)$ is the 
set of all non-crossing partitions of $\{ 1, \ldots , n \}$, and 
``$\leq$'' is reverse refinement order for such partitions. 
$NC(n)$ enters free probability theory via the free cumulants 
introduced in \cite{S94}; a detailed exposition of how this happens 
can be found in Part II of the monograph \cite{NS06}. Both papers 
\cite{BN08}, \cite{D07} build on the basic combinatorics of $NC(n)$, 
but appear to go in different directions, following different goals.

On the one hand, in 2007, Dykema's paper \cite{D07} introduced the 
concept of ``non-crossing linked partition'' and used it to
study the $S$-transform (an important transform introduced by 
Voiculescu \cite{V87} in order to treat the operation of multiplying 
free random variables). In particular, it was shown in \cite{D07} 
how the moment of order $n$ of a noncommutative random variable $a$ 
can be expressed in terms of the coefficients of the series 
$1/S_a (z)$ via a summation over $NCL(n)$, where $S_a (z)$ is the 
$S$-transform of $a$ and where $NCL(n)$ denotes the set of all 
non-crossing linked partitions of $\{ 1, \ldots , n \}$. Following 
to \cite{D07}, further work related to $NCL(n)$ was done in
\cite{CWY08}, \cite{P08}. A review of $NCL(n)$ and of
its connection to the $S$-transform is given in Section 2B and 
Section 4 below.

On the other hand, in 2008, the paper \cite{BN08} by Belinschi and 
Nica introduced a partial order ``$\ll$'' on $NC(n)$, coarser than 
reverse refinement order. The partial order $\ll$ was used in 
\cite{BN08} and subsequently in \cite{BN09} in order to study
connections between free probability and Boolean probability 
(another brand of non-commutative probability, where instead of the 
$NC(n)$'s one works with Boolean posets). A useful fact which 
illustrates how $\ll$ mixes together features from the non-crossing 
and the Boolean worlds is the following: for a fixed $\beta \in NC(n)$, 
the set $\{ \alpha \in NC(n) \mid \alpha \ll \beta \}$ can be 
identified as an interval with respect to reverse refinement order, 
and is consequently counted by a product of Catalan numbers. But 
for a fixed $\alpha \in NC(n)$, the set 
$\{ \beta \in NC(n) \mid \beta \gg \alpha \}$ is 
isomorphic to a Boolean poset, and is thus counted by a power of 2. 
(See review in Propositions \ref{prop:2.4} and \ref{prop:2.5} below.)

The present paper puts into evidence a canonical bijection from
$NCL(n)$ onto the set $\{ ( \alpha , \beta ) \mid  \alpha , 
\beta \in NC(n), \ \alpha \ll \beta \}$. If $\pi \in NCL(n)$ is 
mapped by this bijection to $( \alpha , \beta )$, then:

$\bullet$ $\beta$ is the partition generated by the blocks of $\pi$.
(Note that the blocks of a linked partition are allowed, under 
certain special circumstances, to not be disjoint. Hence $\pi$ does
not have to be itself a partition, and it may happen that
$\beta \neq \pi$.)

$\bullet$ $\alpha$ is obtained by applying a certain ``cycling''
procedure to a partition $\piciup \in NC(n)$ which is defined in 
\cite{D07} and is called the ``unlinking of $\pi$''.

\noindent
The precise details of how this works are given in 
Theorem \ref{thm:3.4} of the paper.

As a consequence of the above bijection, one gets an alternative 
description of Dykema's formula expressing the moments of a 
noncommutative random variable $a$ in terms of the coefficients of 
the reciprocal $S$-transform $1/S_a$. Moreover, due to the Boolean 
features of $\ll$, this formula can be simplified to a form which 
expresses the moment of order $n$ of $a$ by using a summation over 
$NC(n)$ (instead of $NCL(n)$), and where the term indexed by 
$\alpha \in NC(n)$ in the summation can still be canonically 
written as a product over the set of blocks of $\alpha$. The precise 
description of how this goes is given in Theorem \ref{thm:4.5} below. 
It is worth noting that Equation (\ref{eqn:4.51}) from Theorem 
\ref{thm:4.5} has a close resemblance to the moment-cumulant formula 
developped in \cite{BLS96} for the framework of conditionally free 
random variables (another framework where the free and Boolean 
probability worlds interact). 

Besides the introduction, the paper has three other sections. 
Section 2 contains a review of background and notations, 
Section 3 presents the canonical bijection advertised above, and
Section 4 examines the application of this bijection 
to the reciprocal $S$-transform.

$\ $

\begin{center}
{\bf\large 2. Background and notations}
\end{center}
\setcounter{section}{2}
\setcounter{equation}{0}
\setcounter{theorem}{0}

\begin{center}
{\bf 2A. Non-crossing partitions and the partial order
\boldmath{$\ll$} }
\end{center}

\noindent
We will use the standard conventions of notation for non-crossing 
partitions (as in \cite{S00}, or in Lecture 9 of \cite{NS06}). 
The partial order given by reverse refinement on $NC(n)$ will be 
simply denoted as ``$\leq$''; in other words, for 
$\alpha , \beta \in NC(n)$, we write ``$\alpha \leq \beta$'' to 
mean that every block of $\beta$ is a union of blocks of $\alpha$. 
The minimal and maximal element of $( NC(n), \leq )$
are denoted by $0_n$ (the partition of $\{ 1, \ldots , n \}$ into $n$
blocks of 1 element each) and respectively $1_n$ (the partition of
$\{ 1, \ldots , n \}$ into 1 block of $n$ elements).

The coarser partial order $\ll$ on $NC(n)$ was defined in 
\cite{BN08} as follows.

\begin{definition}  \label{def:2.1}
For $\alpha , \beta \in NC(n)$ we write ``$\alpha \ll \beta$'' to
mean that $\alpha \leq \beta$ and that, in addition, the following
condition is fulfilled:
\begin{equation}    \label{eqn:2.11}
\left\{   \begin{array}{l}
\mbox{For every block $W$ of $\beta$ there exists a block}   \\
\mbox{$V$ of $\alpha$ such that $\min (W), \max (W) \in V$.}
\end{array}   \right.
\end{equation}
\end{definition}

It is immediate that if $\alpha \leq \beta$ in $NC(n)$ and if 
$V, W$ are as in (\ref{eqn:2.11}), then one must have 
$V \subseteq W$ and $\min (V) = \min (W)$, $\max (V) = \max (W)$.
It will be convenient to give a special name to the blocks $V$
of $\alpha$ that can be matched to a block $W$ of $\beta$ in 
this way.

\begin{definition}  \label{def:2.2}
Let $\alpha , \beta$ be partitions in $NC(n)$ such that 
$\alpha \ll \beta$. A block $V$ of $\alpha$ will be said to be 
{\em $\beta$-special} when there exists a block $W$ of $\beta$ 
such that $\min (V) = \min (W)$ and $\max (V) = \max (W)$.
\end{definition}

\begin{remark}   \label{rem:2.3}
Let $\alpha , \beta \in NC(n)$ be such that $\alpha \ll \beta$. 

$1^o$ The correspondence $V \mapsto W$ from Definition \ref{def:2.2} 
clearly gives a bijection from the set 
$\{ V \mbox{ block of } \alpha \mid V \mbox{ is $\beta$-special} \}$ 
onto the set of all blocks of $\beta$.

$2^o$ Let us recall that a block $V$ of $\alpha$ is said to be 
{\em inner} (respectively {\em outer}) when there exists (respectively 
when there does not exist) another block $V'$ of $\alpha$ 
such that $\min (V') < \min (V)$ and $\max (V') > \max (V)$.
It is easily seen that every outer block of $\alpha$ is 
$\beta$-special; moreover, the correspondence $V \mapsto W$ from 
Definition \ref{def:2.2} induces a bijection between the outer blocks 
of $\alpha$ and those of $\beta$ -- cf. Remarks 2.9 and 2.12 in 
\cite{BN08}. 
\end{remark}

The next two propositions state in more detail the facts 
mentioned in the introduction about sets of the form 
$\{ \alpha \in NC(n) \mid \alpha \ll \beta \}$ 
and $\{ \beta \in NC(n) \mid \beta \gg \alpha \}$.

\begin{proposition}    \label{prop:2.4}
Let $\beta = \{ W_1, \ldots , W_q \}$ be a partition in $NC(n)$. 
Consider the partition $\beta_0 \in NC(n)$ obtained by refining 
$\beta$ as follows: every block $W_j$ with $|W_j| \leq 2$ is left 
intact, while every block $W_j$ with $|W_j| \geq 3$ is broken into 
the doubleton $\{ \min (W_j), \max (W_j) \}$ and $| W_j | -2$ 
singletons. (Thus every block of $\beta_0$ has either 1 or 2 
elements.) Then 
\begin{equation}   \label{eqn:2.41}
\{ \alpha \in NC(n) \mid \alpha \ll \beta \} = 
\{ \alpha \in NC(n) \mid \beta_0 \leq \alpha \leq \beta \} .
\end{equation}
As a consequence, one has that
\begin{equation}   \label{eqn:2.42}
\vert \ \{ \alpha \in NC(n) \mid \alpha \ll \beta \} \ \vert \ =
\prod_{j=1}^q \cat_{ |W_j|-1 },
\end{equation}
where, for every $k \geq 0$, 
we denote $\cat_k := (2k)!/ \bigl( \, k! (k+1)! \, \bigr)$
(the $k$th Catalan number). 
\end{proposition}

\begin{proof} The equality in (\ref{eqn:2.41}) follows immediately
from how the partial order $\ll$ is defined. The right-hand side of
(\ref{eqn:2.41}) is the interval $[ \beta_0 , \beta ]$ with respect 
to reverse refinement order, and every interval of $(NC(n), \leq )$ 
is known to be canonically isomorphic (as a poset) to a direct 
product of lattices $NC(m)$, $2 \leq m \leq n$; 
see the detailed discussion on pp. 149-150 of
\cite{NS06}, which also gives a concrete algorithm for how 
to obtain the canonical factorization of the interval. By 
following this algorithm it is immediately found that 
\begin{equation}    \label{eqn:2.43}
[ \beta_0, \beta ] \, \simeq \prod_{ \begin{array}{c}
{\scriptstyle 1 \leq j \leq q \ such}   \\
{\scriptstyle that \ |W_j| \geq 3}
\end{array} } \ NC( \, |W_j|-1 \, ),
\end{equation}
and (\ref{eqn:2.42}) follows by taking cardinalities in 
(\ref{eqn:2.43}).
\end{proof}

The next proposition uses the abbreviations
$V \in \alpha$ for ``$V$ is a block of $\alpha$'' and 
$\fV \subseteq \alpha$ for ``$\fV$ is a set of blocks of 
$\alpha$'', where $\alpha$ is a partition in $NC(n)$. 
For the proof of this proposition, the reader is referred 
to Proposition 2.13 and Remark 2.14 of \cite{BN08}. 

\begin{proposition}    \label{prop:2.5}
Let $\alpha$ be in $NC(n)$ and consider the set of partitions 
\begin{equation}  \label{eqn:2.51}
\{ \beta \in NC(n) \mid \beta \gg \alpha \}.
\end{equation}
Then $\beta \mapsto \{ V \in \alpha \mid V
\mbox{ is } \beta{\mbox{-special}} \}$ is a one-to-one map from 
the set (\ref{eqn:2.51}) to the set of subsets of $\alpha$. The 
image of this map is equal to $\{ \fV \subseteq \alpha \mid \fV$
contains all outer blocks of $\alpha \}$.
\end{proposition}

$\ $

\begin{center}
{\bf 2B. Linked partitions and \boldmath{$NCL(n)$} }
\end{center}

\noindent
Following \cite{D07}, we will use the term
``{\em linked partition} of $\{ 1, \ldots , n \}$'' for a set 
$\pi = \{ A_1, \ldots , A_p \}$ 
of non-empty subsets of $\{ 1, \ldots , n \}$ such that 
$A_1 \cup \cdots \cup A_p = \{ 1, \ldots , n \}$, and where 
for every $i \neq j$ ($1 \leq i,j \leq p$) one has that either 
$A_i \cap A_j = \emptyset$ or that the following holds:
\begin{equation}    \label{eqn:2B1}
\left\{   \begin{array}{l}
| A_i | \geq 2, \, | A_j | \geq 2, \, | A_i \cap A_j | =1, 
\ \min (A_i) \neq \min (A_j), \mbox{ and }  \\
\mbox{the unique element of 
$A_i \cap A_j$ is one of $\min (A_i)$, $\min (A_j)$.}
\end{array}  \right.
\end{equation}
Obviously, every partition of $\{ 1, \ldots , n \}$ (in the usual 
sense of the term) is a linked partition of $\{ 1, \ldots , n \}$, 
but the converse is not true. Throughout this paper we use the 
letters $\alpha , \beta , \ldots$ to denote partitions, and 
$\pi, \rho, \ldots$ to denote linked partitions (which may or 
may not be partitions). A few more terms and basic facts 
from Section 5 of \cite{D07} are reviewed next.

\begin{review}   \label{rev:2.6}
$1^o$ If $\pi = \{ A_1, \ldots , A_p \}$ is a linked partition 
of $\{ 1, \ldots , n \}$ then $A_1, \ldots , A_p$ are called the 
{\em blocks} of $\pi$. It is easy to see that
every $m \in \{ 1, \ldots , n \}$ belongs either to exactly one or 
to exactly two blocks $A_i$; in the first case one says that $m$ is 
{\em singly-covered} by $\pi$, and in the second case one says that 
$m$ is {\em doubly-covered} by $\pi$.

$2^o$ Let $\pi$ be a linked partition of $\{ 1, \ldots , n \}$. 
The partition of $\{ 1, \ldots , n \}$ which is generated by $\pi$
will be denoted as $\widehat{\pi}$. In other words, 
$\widehat{\pi}$ is the smallest (with respect to reverse refinement
order) among all partitions $\beta$ of $\{ 1, \ldots , n \}$ which
have the following property: ``for every block $A$ of $\pi$ there 
exists a block $V$ of $\beta$ such that $A \subseteq V$''.

$3^o$ Let $\pi = \{ A_1, \ldots , A_p \}$ be a linked partition of 
$\{ 1, \ldots , n \}$. For every $1 \leq i \leq p$ define 
\begin{equation}
V_i = \left\{  \begin{array}{ll}
A_i,  & \mbox{if $\min (A_i)$ is singly-covered by $\pi$}   \\
A_i \setminus \{ \min (A_i),  & \mbox{if $\min (A_i)$ is 
                                 doubly-covered by $\pi$.}   
\end{array}   \right.
\end{equation}
Then $\{ V_1, \ldots , V_p \}$ is a partition of $\{1, \ldots , n \}$,
called the {\em unlinking} of $\pi$ and denoted as $\piciup$.

$4^o$ A linked partition $\pi$ of $\{ 1, \ldots , n \}$ is said to 
be {\em non-crossing} if it is not possible to find two distinct 
blocks $A,B$ of $\pi$ and elements $a,a' \in A$, $b,b' \in B$ such 
that $a<a'<b<b'$. 

$5^o$ Let $n$ be a positive integer. The set of all non-crossing 
linked partitions of $\{ 1, \ldots , n \}$ is denoted by $NCL(n)$.
It is not hard to see that if $\pi \in NCL(n)$ then the partitions
$\widehat{\pi}$ and $\piciup$ defined in $2^o$ and $3^o$ above 
belong to $NC(n)$.
\end{review}

\begin{remark}    \label{rem:2.7}
{\em (Restrictions of linked partitions.)}
It is immediate that the above discussion about linked partitions
could be carried out without any modifications in the larger 
framework where instead of $\{ 1, \ldots , n \}$ one uses an
abstract finite totally ordered set $F$, and one considers the
set $NCL(F)$ of non-crossing linked partitions of $F$ (instead 
of just sticking to the $NCL(n)$). This doesn't really bring 
anything new, since it is obvious that $NCL(F)$ is canonically 
identified to $NCL( \, |F| \, )$ via the map 
$\{ A_1, \ldots , A_p \} \mapsto \{ f(A_1), \ldots , f(A_p) \}$,
where $f$ is the unique order-preserving bijection from $F$ onto 
$\{ 1, \ldots , |F| \}$. But in the subsequent discussion it 
will be nevertheless convenient to allow linked partitions for
a slightly more general kind of set $F$ (specifically, for $F$ 
being a non-empty finite subset of $\bN$), in order to 
simplify the notations for restrictions of linked partitions.

So let $E \subseteq F$ be non-empty subsets of $\bN$, and let 
$\pi \in NCL(F)$ be such that $E$ is saturated with respect to 
$\pi$ (which means that whenever $A$ is a block of $\pi$ and
$A \cap E \neq \emptyset$, it follows that $A \subseteq E$). 
Then $\pi$ is of the form 
$\{ A_1, \ldots , A_p, A_1', \ldots , A_q' \}$
with $A_1, \ldots , A_p \subseteq E$ and
$A_1', \ldots , A_q' \subseteq F \setminus E$, and one defines 
the {\em restriction} of $\pi$ to $E$ to be
\begin{equation}   \label{eqn:2.71}
\pi \, \vert \, { }_E := \{ A_1, \ldots , A_p \} .
\end{equation}
It is immediate that $\pi \, \vert \, { }_E \in NCL(E)$. Moreover, 
it is easily seen that the operation of restriction is 
well-behaved with respect to the maps $\pi \mapsto \widehat{\pi}$ 
and $\pi \mapsto \piciup$ from $2^o$ and $3^o$ of Review 
\ref{rev:2.6}, in the sense that one has
\begin{equation}   \label{eqn:2.72}
\Bigl( \, \pi \, \vert \, { }_E \, \Bigr)^{\wedge}  = 
\bigl( \, \widehat{\pi} \, \bigr) \, \vert \, { }_E , 
\ \ \mbox{ and } \ \ 
\Bigl( \, \pi \, \vert \, { }_E \, \Bigr)^{\vee}  = 
\bigl( \, \piciup \, \bigr) \, \vert \, { }_E  .
\end{equation}
\end{remark}

Let us next record the observation (cf. \cite{D07}, 
Corollary 5.13 and its proof) that non-crossing linked partitions 
can be broken into ``irreducible'' pieces, as follows.

\begin{proposition}  \label{prop:2.8}
Let $n$ be a positive integer and let 
$\beta = \{ W_1, \ldots , W_q \}$ be a partition in $NC(n)$.
The map
\[
\pi \mapsto \Bigl( \pi \, \vert \, { }_{W_1}, \ldots , 
\pi \, \vert \, { }_{W_q} \Bigr)
\]
is a bijection from 
$\{ \pi \in NCL(n) \mid \widehat{\pi} = \beta \}$ onto
$\prod_{j=1}^q
\{ \pi_j \in NCL(W_j) \mid \widehat{\pi}_j = 1_{W_j} \}$
(where $1_{W_j} \in NC( W_j )$ is the partition of $W_j$ into 
only one block).
\end{proposition}

So, after suitably renumbering every $W_j$ from 
the preceding proposition as $\{ 1, \ldots , |W_j| \}$, 
one is reduced in the end to looking at sets of 
non-crossing linked partitions of the form 
$\{ \pi \in NCL(m) \mid \widehat{\pi} = 1_m \}$.
In Proposition 5.11 of \cite{D07} it is pointed out that 
the latter sets can be identified with sets of usual 
non-crossing partitions, as follows.

\begin{proposition}  \label{prop:2.9}
For every $n \geq 2$, the map $\pi \mapsto \piciup$ 
is a bijection from 
$\{ \pi \in NCL(n) \mid \widehat{\pi} = 1_n \}$ onto
$\{ \alpha \in NC(n) \mid 1 \mbox{ and } 2$
belong to the same block of $\alpha \}$.
\end{proposition}

$\ $

\begin{center}
{\bf\large 3. The canonical bijection relating \boldmath{$NCL(n)$} 
to \boldmath{$\ll$} }
\end{center}
\setcounter{section}{3}
\setcounter{equation}{0}
\setcounter{theorem}{0}

\noindent
This section is devoted to proving the bijection announced
in the introduction. The main result is Theorem \ref{thm:3.4}.
We start by recalling some connections between set-partitions 
and permutations, and by defining precisely what is the 
``cycled unlinking'' of a linked partition $\pi \in NCL(n)$.

\begin{notation}    \label{def:3.1}
$1^o$ If $\tau$ is a permutation of $\{ 1, \ldots , n \}$ and 
if $\alpha = \{ V_1, \ldots , V_p \}$ is a partition of 
$\{ 1, \ldots , n \}$ then we denote
\begin{equation}   \label{eqn:3.11}
\tau \cdot \alpha := \{ \tau (V_1), \ldots , \tau (V_p) \}
\ \ \mbox{ (a new partition of $\{ 1, \ldots , n \}$). }
\end{equation}
It is clear that formula (\ref{eqn:3.11}) defines an action
of the symmetric group $\cS_n$ on the set of all partitions of 
$\{ 1, \ldots , n \}$.

$2^o$ Every partition $\alpha \in NC(n)$ has associated to it a
permutation of $\{ 1, \ldots , n \}$, which is denoted by
$P_{\alpha}$, and is defined by the following prescription: for
every block $V = \{ i_1, \ldots , i_m \}$ of $\pi$, with
$i_1 < \cdots < i_m$, one creates a cycle of $P_{\alpha}$ by 
putting
\[
P_{\alpha} (i_1) = i_2, \ldots , P_{\alpha} (i_{m-1}) = i_m,
P_{\alpha} (i_m) = i_1.
\]
\end{notation}

\begin{definition}    \label{def:3.2}
Let $\pi$ be in $NCL(n)$. Consider the partitions 
$\widehat{\pi}, \piciup \in NC(n)$ (as in Review \ref{rev:2.6}),
and form the new partition 
\begin{equation}    \label{eqn:3.21}
\picirc \ := \ P_{\widehat{\pi}}^{-1} \cdot \piciup 
\end{equation}
(where the permutation $P_{\widehat{\pi}}$ is defined as in 
Notation \ref{def:3.1}.2). The partition $\picirc$ will be called 
the {\em cycled unlinking} of $\pi$.
\end{definition}

\begin{example}   \label{ex:3.3}
{\em (A concrete example.)}
Say for instance that $n=11$ and that 
\[
\pi = \bigl\{ \, \{ 1,2,4 \} , \{ 2,3 \} , \{ 4,5,6 \},
\{ 6,7 \} , \{ 8,9,11 \}, \{ 9,10 \} \, \bigr\} \in NCL(11).
\]
Then one has
\[
\widehat{\pi} = 
\bigl\{ \, \{ 1,2,3,4,5,6,7 \} , \{ 8,9,10,11 \} \, \bigr\} , 
\ \ \piciup = \bigl\{ \, \{ 1,2,4 \} , \{ 3 \} , \{ 5,6 \},
\{ 7 \} , \{ 8,9,11 \}, \{ 10 \} \, \bigr\} .
\]
The permutation associated to $\widehat{\pi}$ is 
\[
P_{\widehat{\pi}} = \bigl( 1,2,3,4,5,6,7 \bigr) \
\bigl( 8,9,10,11 \bigr) \ \ \mbox{ (written in cycle notation),}
\]
hence the cycled unlinking of $\pi$ is 
\[
\picirc \, = P_{\widehat{\pi}}^{-1} \cdot \piciup \, =
\bigl\{ \, \{ 1,3,7 \} , \{ 2 \} , \{ 4,5 \},
\{ 6 \} , \{ 8,10,11 \}, \{ 9 \} \, \bigr\} .
\]
Note that $\picirc \, \in NC(11)$ and $\picirc \, \ll \widehat{\pi}$
(as we will see that it must generally be the case).
\end{example}

\begin{theorem}  \label{thm:3.4}
Let $n$ be a positive integer. The map 
$\pi \mapsto \bigl( \, \picirc , \widehat{\pi} \, \bigr)$
is a bijection from $NCL(n)$ onto the set
$\{ ( \alpha , \beta ) \mid \alpha , \beta \in NC(n), \
\alpha \ll \beta \}$.
\end{theorem}

For the proof of Theorem \ref{thm:3.4}, it is convenient to 
first establish some lemmas.

\begin{lemma}   \label{lemma:3.5}
Let $\beta$ be a partition in $NC(n)$. The map 
$\alpha \mapsto P_{\beta}^{-1} \cdot \alpha$ sends the set
$\{ \alpha \in NC(n) \mid \alpha \leq \beta \}$ into itself.
\end{lemma}

\begin{proof} Let us fix an $\alpha \in NC(n)$ such that 
$\alpha \leq \beta$, and let us consider the partition
$\alpha ' := P_{\beta}^{-1} \cdot \alpha$. Observe that
$\alpha ' \leq \beta$ in reverse refinement order, i.e that every 
block $V'$ of $\alpha '$ is contained in a block of $\beta$. Indeed, 
$V'$ can be written as $V' = P_{\beta}^{-1} (V)$ where $V$ is a 
block of $\alpha$, and the block $W$ of $\beta$ which contains $V$ 
must contain $V'$ as well (since 
$V' \subseteq P_{\beta}^{-1} (W) = W$).

It remains to show that the 
partition $\alpha '$ is non-crossing. Let us fix two distinct 
blocks $V_1 '$ and $V_2 '$ of $\alpha '$, and let us prove that 
$V_1'$ and $V_2'$ do not cross. Consider the blocks $W_1, W_2$
of $\beta$ such that $W_1 \supseteq V_1'$ and $W_2 \supseteq V_2'$, 
and distinguish the following two cases.

\noindent
{\em Case 1:} $W_1 \neq W_2$. 
In this case $W_1$ and $W_2$ do not cross, so the smaller sets
$V_1'$ and $V_2'$ can't cross either. 

\noindent
{\em Case 2:} $W_1 = W_2 =:W$.
In this case, consider the sets $V_1 := P_{\beta} (V_1')$ and 
$V_2 := P_{\beta} (V_2')$. These are two distinct blocks of $\alpha$, 
which are still contained in $W$. We know that $V_1$ and $V_2$ don't 
cross (since $\alpha$ is in $NC(n)$), and we will conclude the proof
by arguing that if $V_1'$ and $V_2'$ would cross, then so would 
$V_1$ and $V_2$. Indeed, suppose that $i_1, j_1 \in V_1'$ and 
$i_2, j_2 \in V_2'$ are such that $i_1 < i_2 < j_1 < j_2$. If $j_2$ 
is not the maximal element of $W$, then it follows that 
$P_{\beta} (i_1) < P_{\beta} (i_2) < P_{\beta} (j_1) 
< P_{\beta} (j_2)$, 
which is a crossing between $V_1$ and $V_2$. In the opposite 
situation when $j_2 = \max (W)$, it follows that 
$P_{\beta} (j_2) < P_{\beta} (i_1) < 
P_{\beta} (i_2) < P_{\beta} (j_1)$, 
and we get a crossing between $V_2$ and $V_1$ (since
$P_{\beta} (j_2) , P_{\beta} (i_2)  \in V_2$ and
$P_{\beta} (i_1) , P_{\beta} (j_1)  \in V_1$). 
\end{proof}

\begin{lemma}   \label{lemma:3.6}
Let $\pi$ be in $NCL(n)$. Then the partition 
$\picirc$ is in $NC(n)$, and 
$\picirc \ll \widehat{\pi}$.
\end{lemma}

\begin{proof} Since $\piciup, \widehat{\pi} \in NC(n)$ and 
$\piciup \leq \widehat{\pi}$, the preceding lemma gives us that 
$\picirc = P_{\widehat{\pi}}^{-1} \cdot \piciup$ is in $NC(n)$ and 
that $\picirc \leq \widehat{\pi}$. 

It remains to show that, for every block $W$ of $\widehat{\pi}$,
the numbers $\min (W)$ and $\max (W)$ belong to the same block of
$\picirc$. For the rest of the proof we fix such a $W$, with 
$|W| =: m \geq 2$ (the case $|W| =1$ is obvious). Let us write 
explicitly $W = \{ i_1, i_2, \ldots , i_m \}$, with 
$i_1 < i_2 < \cdots < i_m$; our goal is then to prove that $i_1$ 
and $i_m$ belong to the same block of $\picirc$.

Consider the linked partition $\pi \ \vert \ { }_W \in NCL(W)$, 
which has 
$\bigl( \, \pi \ \vert \ { }_W \, \bigr)^{\wedge}$ = 
$\widehat{\pi} \ \vert \ { }_W = 1_W$. Proposition \ref{prop:2.9}
gives us (after going through the suitable identification of 
$i_1, i_2, \ldots , i_m$ with $1,2, \ldots , m$) that $i_1$ and 
$i_2$ belong to the same block of 
$\bigl( \, \pi \ \vert \ { }_W \, \bigr)^{\vee}$.
But $\bigl( \, \pi \ \vert \ { }_W \, \bigr)^{\vee}$ =
$\piciup \ \vert \ { }_W$ (cf. Equation (\ref{eqn:2.72})), hence
it follows that $i_1$ and $i_2$ belong to the same block $V$ of 
$\piciup$. This implies in turn that $P_{\widehat{\pi}}^{-1} (i_1)$ 
and $P_{\widehat{\pi}}^{-1} (i_2)$ belong to the same block 
$P_{\widehat{\pi}}^{-1} (V)$ of $\picirc$; but  
$P_{\widehat{\pi}}^{-1} (i_1) = i_m$ and 
$P_{\widehat{\pi}}^{-1} (i_2) = i_1$, so we are done. 
\end{proof}

\begin{lemma}   \label{lemma:3.7}
Let $\pi, \rho \in NCL(n)$ be such that 
$\widehat{\pi} = \widehat{\rho}$ and 
$\piciup = \rhociup$. Then $\pi = \rho$.
\end{lemma}

\begin{proof} Let us denote 
$\widehat{\pi} = \widehat{\rho} =: \beta \in NC(n)$, and let us
write explicitly $\beta = \{ W_1, \ldots , W_q \}$. In order to 
prove that $\pi = \rho$ it suffices to verify that they have the 
same image by the bijection from Proposition \ref{prop:2.8}, i.e.
that $\pi \ \vert \ { }_{W_j} = \rho \ \vert \ { }_{W_j} 
\in NCL(W_j)$ for every $1 \leq j \leq q$.

So let us fix $j$ ($1 \leq j \leq q$). Observe that 
$\bigl( \, \pi \ \vert \ { }_{W_j} \, \bigr)^{\wedge}$
= $\widehat{\pi} \ \vert \ { }_{W_j}$ = 
$\beta \ \vert \ { }_{W_j}$ = $1_W$, and similarly
$\bigl( \, \rho \ \vert \ { }_{W_j} \, \bigr)^{\wedge}$
= $1_{W_j}$. Now, by invoking Proposition \ref{prop:2.9} (and by 
using the suitable renumbering of $W_j$ into
$\{ 1, \ldots , |W_j| \}$) we see that the map 
$\theta \mapsto \thetaciup$ is one-to-one on 
$\{ \theta \in NCL(W_j) \mid \widehat{\theta} = 1_{W_j} \}$.
Thus the required fact that 
$\pi \ \vert \ { }_{W_j} = \rho \ \vert \ { }_{W_j}$
will follow if we can prove that 
$\bigl( \, \pi \ \vert \ { }_{W_j} \, \bigr)^{\vee}$
= $\bigl( \, \rho \ \vert \ { }_{W_j} \, \bigr)^{\vee}$.
But the latter equality amounts (in view of (\ref{eqn:2.72})) to 
$\piciup \ \vert \ { }_{W_j} = \rhociup \ \vert \ { }_{W_j}$,
and thus follows from the hypothesis that $\piciup = \rhociup$.
\end{proof}

\vspace{10pt}

\noindent
{\em Proof of Theorem \ref{thm:3.4}.}
From Lemma \ref{lemma:3.6} it follows that the map 
$\pi \mapsto ( \picirc , \widehat{\pi} )$ 
is well-defined (with target set as described in the theorem).
In order to prove the injectivity of this map, consider two 
linked partitions $\pi , \rho \in NCL(n)$ such that 
$\picirc \, = \, \rhocirc \, =: \alpha$ and 
$\widehat{\pi} = \widehat{\rho} =: \beta$. Then 
$\piciup \, = \rhociup \, = P_{\beta} \cdot \alpha$; hence 
$\pi$ and $\rho$ satisfy the hypotheses of Lemma \ref{lemma:3.7},
and it follows that $\pi = \rho$.

In order to complete the proof of the theorem, it now suffices 
to verify that the sets $NCL(n)$ and 
$\{ ( \alpha , \beta ) \mid \alpha , \beta \in NC(n), \
\alpha \ll \beta \}$
have the same cardinality. For $NCL(n)$ we can invoke Corollary 5.13 
of \cite{D07}, which gives us the formula
\begin{equation}    \label{eqn:3.41}
| \, NCL(n) \, | \ = \ \sum_{\beta \in NC(n)} \
\Bigl( \, \prod_{W \in \beta}  \cat_{ |W|-1 } \, \Bigr) . 
\end{equation}
But on the other hand we have
\[
\vert \ \{ ( \alpha , \beta ) \mid \alpha , \beta \in NC(n), \
\alpha \ll \beta \} \ \vert \ 
= \sum_{\beta \in NC(n)} \ \vert \ \{ \alpha 
\mid \alpha \in NC(n), \ \alpha \ll \beta \} \ \vert \ ,
\]
which is indeed equal to the right-hand side 
of (\ref{eqn:3.41}), due to Proposition \ref{prop:2.4}.
\hfill{$\blacksquare$}

\vspace{6pt}

\begin{remark}   \label{rem:3.8}
Theorem \ref{thm:3.4} allows one to transfer enumeration properties
in between the two sets involved in the bijection of the theorem.
In particular, as it is known that $NCL(n)$ is counted by the 
$(n-1)$th Schr\"oder number $r_{n-1}$ (Theorem 8.3 in \cite{D07}, 
see also the discussion in Section 2 of \cite{CWY08}), it follows 
that the same is true for $\{ ( \alpha , \beta ) \mid$
$\alpha , \beta \in NC(n), \ \alpha \ll \beta \}$. By writing the 
latter set as 
$\cup_{\alpha} \{ \beta \in NC(n) \mid \beta \gg \alpha \}$
and by invoking Proposition \ref{prop:2.5} one thus comes to the
following amusing (possibly already known) enumerative interpretation 
of the Schr\"oder number $r_{n-1}$: it counts coloured 
non-crossing partitions $\alpha \in NC(n)$, where every block of 
$\alpha$ is coloured in either red or blue, and all the outer blocks 
are red. For a nice collection of other enumerative interpretations 
of the Schr\"oder numbers, see Example 6.2.8 and Exercise 6.39 of 
the monograph \cite{S99}.
\end{remark}

$\ $

\begin{center}
{\bf\large 4. Application to the reciprocal 
\boldmath{$S$}-transform}
\end{center}
\setcounter{section}{4}
\setcounter{equation}{0}
\setcounter{theorem}{0}

\noindent
Let $( \cA , \varphi )$ be an algebraic noncommutative probability 
space (which simply means that $\cA$ is a unital algebra over $\bC$ 
and $\varphi : \cA \to \bC$ is a linear functional such that 
$\varphi ( 1_{\cA} ) = 1$), and let $a \in \cA$ be such that 
$\varphi (a) = 1$. The {\em $S$-transform} of $a$ is the series
\begin{equation}    \label{eqn:4.1}
S_a (z) := \frac{1+z}{z} M_a^{\langle -1 \rangle} (z)
\in \bC \, [[ z ]],
\end{equation}
where $M_a (z) := \sum_{n=1}^{\infty} \varphi (a^n) z^n$ (the 
{\em moment series} of $a$), and where $M_a^{\langle -1 \rangle}$
is the inverse of $M_a$ under composition. The $S$-transform plays
a fundamental role in the study of multiplication of free random 
variables (see Section 3.6 of \cite{VDN92} or Lecture 18 of 
\cite{NS06} for the details of how this goes).

It is immediate that the series $S_a (z)$ from (\ref{eqn:4.1})
has constant term equal to 1. Hence one can consider its reciprocal
(i.e. inverse under multiplication) $1/S_a (z)$, which is another 
series with constant term equal to 1:
\begin{equation}  \label{eqn:4.2}
1/S_a (z) = \sum_{n=0}^{\infty} t_n z^n, \ \ \ \mbox{ with }
t_0 = 1.
\end{equation}
In the paper \cite{D07} the series $1/S_a$ goes under the name of 
``{\em $T$-transform of $a$}''. It is observed there that the formula
giving back the family of moments $\{ \varphi (a^n) \mid n \geq 1 \}$
of $a$ in terms of the family $\{ t_n \mid n \geq 0 \}$ of 
coefficients of $1/S_a$ only uses positive integer coefficients:
\begin{equation}   \label{eqn:4.3}
\left\{   \begin{array}{l}
\varphi (a) = 1, \ \varphi (a^2) = t_1 + 1, \
\varphi (a^3) = t_2 + t_1^2 + 3t_1 +1,               \\
\varphi (a^4) = t_3 + 3t_2 t_1 + t_1^3 + 4t_2 + 6t_1^2 + 6t_1 + 1,
, \ldots 
\end{array}  \right.
\end{equation}
Moreover, \cite{D07} identifies precisely the combinatorial structure
which indexes the sums in (\ref{eqn:4.3}): this is nothing but
$NCL(n)$, and the formula generalizing the special cases from 
(\ref{eqn:4.3}) is
\begin{equation}   \label{eqn:4.4}
\varphi (a^n) = \sum_{\pi \in NCL(n)} \ \Bigl( \, 
\prod_{ A \in \pi } \, t_{ |A|-1 } \, \Bigr) , \ \ n \geq 1,
\end{equation}
where ``$A \in \pi$'' is an abbreviation for ``$A$ is a block 
of $\pi$''; see Proposition 8.1 in \cite{D07}.

In view of Theorem \ref{thm:3.4} of the present paper, the formula
giving $\varphi (a^n)$ in terms of the coefficients of $1/S_a$ can 
be equivalently understood as a summation over the set 
$\{ ( \alpha , \beta ) \mid \alpha , \beta \in NC(n)$,
$\alpha \ll \beta \}$. This alternative version of the formula will 
be stated in Proposition \ref{prop:4.3} below. In preparation of 
that statement, we prove two lemmas.

\begin{lemma}   \label{lemma:4.1}
Let $\pi$ be a linked partition in $NCL(n)$ such that 
$\widehat{\pi} = 1_n$ and let $A$ be a block of $\pi$ such that 
$A \not\ni 1$. Then $\min (A)$ is doubly-covered by $\pi$.
\end{lemma}

\begin{proof} Denote $\min (A) =: m$. Since $1$ and $m$ belong to 
the same (unique) block of $\widehat{\pi}$, there have to exist 
$p \geq 1$, some 
$m_0, m_1, \ldots , m_p \in \{ 1, \ldots , n \}$ and some blocks 
$A_1, \ldots , A_p$ of $\pi$ such that $m_0 =1$, $m_p =m$, and 
such that
\begin{equation}   \label{eqn:4.11}
m_0, m_1 \in A_1; \, m_1, m_2 \in A_2; \, \ldots , 
m_{p-1}, m_p \in A_p.
\end{equation}
Let us suppose moreover that in (\ref{eqn:4.11}) $p$ is picked 
to be as small as possible. It is then immediate that 
$m_{j-1} \neq m_j$ for every $1 \leq j \leq p$, and that 
$A_{j-1} \neq A_j$ for every $2 \leq j \leq p$. 

Observe that $m_1 \in A_1 \cap A_2$ and $m_1$ is not the minimum 
of $A_1$ (since $m_0 \in A_1$, and $m_0 = 1 < m_1$); from the 
definition of a linked partition it follows that $m_1 = \min (A_2)$.
We next observe that $m_2 \in A_2 \cap A_3$ and 
$m_2 \neq \min (A_2)$ (since $\min (A_2) = m_1 \neq m_2$), so the 
same argument as above applies to give us that $m_2 = \min (A_3)$. 
Continuing like this by induction we find that 
$m_j = \min (A_{j+1})$ for every $1 \leq j \leq p-1$, and in 
particular that $m_{p-1} = \min (A_p)$.
\footnote{ This argument was run by assumming that $p \geq 2$. If 
$p=1$ then the conclusion that $m_{p-1} = \min (A_p)$ still holds,
as we must have that $m_0 = 1 = \min (A_1)$. }
Thus the block $A_p$ of $\pi$ contains $m_p = m$, and $A_p$ is 
different from $A$ (because 
$\min (A_p) = m_{p-1} \neq m_p = \min (A)$). This shows that $m$ is
doubly-covered by $\pi$, as required.
\end{proof}

\begin{lemma}   \label{lemma:4.2}
Let $\pi$ be a linked partition in $NCL(n)$, let $A$ be a block 
of $\pi$, and let $W$ be the unique block of $\widehat{\pi}$
such that $W \supseteq A$. Then
\begin{equation}    \label{eqn:4.21}
\left(  \begin{array}{c}
\mbox{$\min (A)$ is}       \\
\mbox{singly-covered by $\pi$}
\end{array}   \right) \ \Leftrightarrow \
\Bigl( \, \min (A) = \min (W) \, \Bigr) .
\end{equation}
\end{lemma}

\begin{proof}
By replacing $\pi$ with $\pi \, \vert \, { }_W$ and by redenoting 
the elements of $W$ as $\{ 1, 2 , \ldots , |W| \}$ in increasing 
order, we may assume without loss of generality that 
$\widehat{\pi} = 1_n$ and hence that $W = \{ 1, \ldots , n \}$. 
The statement on the right-hand side of equivalence (\ref{eqn:4.21}) 
becomes ``$\min (A) = 1$''. The implication ``$\Leftarrow$'' in this 
equivalence is then immediate (as the definition of a 
linked partition implies that $1$ always is singly-covered), while 
the implication ``$\Rightarrow$'' follows from Lemma \ref{lemma:4.1}.
\end{proof}

\begin{proposition}   \label{prop:4.3}
Let $( \cA , \varphi )$ be an algebraic noncommutative probability 
space, let $a \in \cA$ be such that $\varphi (a) = 1$, and consider 
the reciprocal $S$-transform 
$1/S_a (z) = \sum_{n=0}^{\infty} t_n z^n$. Then for every $n \geq 1$
one has
\begin{equation}   \label{eqn:4.31}
\varphi (a^n) = \sum_{ \begin{array}{c}
{\scriptstyle \alpha , \beta \in NC(n)}   \\
{\scriptstyle such \ that \ \alpha \ll \beta}
\end{array} } \ 
\Bigl( \, \prod_{ \begin{array}{c}
{\scriptstyle U \in \alpha,} \\
{\scriptstyle \beta-special}
\end{array} } \ t_{ |U|-1 } \, \Bigr) \cdot 
\Bigl( \, \prod_{ \begin{array}{c}
{\scriptstyle V \in \alpha,} \\
{\scriptstyle not \ \beta-special}
\end{array} } \ t_{ |V| } \, \Bigr) 
\end{equation}
(where the concept of $\beta$-special
block of $\alpha$ is as in Definition \ref{def:2.2}).
\end{proposition}

\begin{proof} 
We will verify that the sums on the right-hand sides of 
(\ref{eqn:4.4}) and (\ref{eqn:4.31}) are identified term by term
when one uses the bijection $\pi \leftrightarrow ( \alpha , \beta )$
from Theorem \ref{thm:3.4}. We thus fix for the whole proof 
$\pi \in NCL(n)$ and $\alpha , \beta \in NC(n)$ such that 
$\pi \leftrightarrow ( \alpha , \beta )$ in Theorem \ref{thm:3.4},
and our goal is to show that 
\begin{equation}   \label{eqn:4.32}
\prod_{ A \in \pi } \ t_{ |A|-1 } \ = \
\Bigl( \, \prod_{ \begin{array}{c}
{\scriptstyle U \in \alpha,} \\
{\scriptstyle \beta-special}
\end{array} } \ t_{ |U|-1 } \, \Bigr) \cdot 
\Bigl( \, \prod_{ \begin{array}{c}
{\scriptstyle V \in \alpha,} \\
{\scriptstyle not \ \beta-special}
\end{array} } \ t_{ |V| } \, \Bigr) .
\end{equation}

The fact that $\pi \leftrightarrow ( \alpha , \beta )$ in Theorem 
\ref{thm:3.4} means of course that $\alpha = \, \picirc$ and
$\beta = \widehat{\pi}$. Let us write explicitly 
$\beta = \{ W_1, \ldots , W_q \}$. In view of Lemma \ref{lemma:4.2} 
we see that $\pi$ can then be written in the form 
$\pi = \{ A_1, \ldots , A_q, B_1, \ldots , B_r \}$
where $\min (A_1) = \min (W_1), \ldots , \min (A_q) = \min (W_q)$
and where $\min (A_1), \ldots , \min (A_q)$ are singly-covered by 
$\pi$, while $\min (B_1), \ldots , \min (B_r)$ are doubly-covered 
by $\pi$. From how the unlinking $\piciup$ is defined (cf.
Review \ref{rev:2.6}.3) we next infer that 
\[
\left\{   \begin{array}{l}
\piciup = \{ U_1, \ldots , U_q, V_1, \ldots , V_r \} , 
\ \ \mbox{ where }                                       \\
U_j = A_j, \ \forall \, 1 \leq j \leq q \mbox{ and }
V_k = B_k \setminus \{ \min (B_k) \}, \ \forall \, 1 \leq k \leq r .
\end{array}   \right.
\]
From Definition \ref{def:3.2} we further infer that the 
partition $\alpha = \, \picirc$ can be written in the form
\[
\left\{   \begin{array}{l}
\picirc = \{ U_1', \ldots , U_q', V_1', \ldots , V_r' \} , 
\ \ \mbox{ where }                                       \\
U_j' = P_{\beta}^{-1}(U_j), \ \forall \, 1 \leq j \leq q \mbox{ and }
V_k' = P_{\beta}^{-1}(V_k), \ \forall \, 1 \leq k \leq r .
\end{array}   \right.
\]

We next observe that 
\begin{equation}   \label{eqn:4.33}
\min (W_j) , \max (W_j) \in U_j' , \ \ \forall \, 1 \leq j \leq q.
\end{equation}
Indeed, this statement is clear in the case when $|W_j| = 1$ (and
when $W_j = U_j = U_j'$). In the case when $|W_j| =:m \geq 2$
we write $W_j = \{ i_1, i_2, \ldots , i_m \}$ with 
$i_1 < i_2 < \cdots < i_m$ and we observe that the argument used in 
the proof of Lemma \ref{lemma:3.6} applies, giving us that 
$i_1, i_2 \in U_j$, and hence that $i_1, i_m \in U_j'$. 

From (\ref{eqn:4.33}) it follows that the $\beta$-special blocks of 
$\alpha$ are precisely $U_1', \ldots , U_q'$. The right-hand side 
of (\ref{eqn:4.32}) thus takes the form 
\begin{equation}   \label{eqn:4.34}
\Bigl( \, \prod_{j=1}^q \ t_{ |U_j'|-1 } \, \Bigr) \cdot 
\Bigl( \, \prod_{k=1}^r \ t_{ |V_k'| } \, \Bigr) .
\end{equation}
But it is clear that $|U_j'| = |U_j| = |A_j|$, 
$1 \leq j \leq q$, and that $|V_k'| = |V_k| = |B_k| - 1$, 
$1 \leq k \leq r$.  Hence the product from (\ref{eqn:4.34}) equals
$\bigl( \, \prod_{j=1}^q \ t_{ |A_j|-1 } \, \bigr) \cdot 
\bigl( \, \prod_{k=1}^r \ t_{ |B_k|-1 } \, \bigr)$,
and (\ref{eqn:4.32}) follows.
\end{proof}

It is in fact fairly easy, in hindsight, to give a direct proof 
of Proposition \ref{prop:4.3} by using the $R$-transform $R_a$
(another important transform of free probability) as an intermediate 
for passing from the series $1/S_a$ to the moments of $a$.

\vspace{14pt}

\noindent
{\em Second proof of Proposition \ref{prop:4.3}.}
Let us consider the $R$-transform of $a$. This is the series
$R_a (z) = \sum_{n=1}^{\infty} \kk_n z^n$ whose coefficients 
$( \kk_n )_{n=1}^{\infty}$, called the {\em free cumulants} of $a$,
are uniquely determined by the fact that they satisfy the relations
\begin{equation}  \label{eqn:4.35}
\varphi (a^n) = \sum_{\beta \in NC(n)} \
\Bigl( \, \prod_{W \in \beta} \, \kk_{|W|} \, \Bigr) , \ \ 
\forall \, n \geq 1 .
\end{equation}
(For instance $\kk_1 = \varphi (a) \bigl( =1 \bigr)$,
$\kk_2 = \varphi (a^2) - \bigr( \varphi (a) \bigr)^2$ and
$\kk_3 = \varphi (a^3) - 3 \varphi (a) \, \varphi (a^2) +
2 \bigr( \varphi (a) \bigr)^3$. For a more detailed discussion of
Equation (\ref{eqn:4.35}), see e.g. Lecture 12 of \cite{NS06}.)
It was recently observed in \cite{MN08} that the formula expressing
the free cumulants of $a$ in terms of the coefficients of $1/S_a (z)$ 
is merely a shifted version of Equation (\ref{eqn:4.35}):
\begin{equation}  \label{eqn:4.36}
\kk_n = \sum_{\gamma \in NC(n-1)} \
\Bigl( \, \prod_{V \in \gamma} \, t_{|V|} \, \Bigr) , \ \ 
\forall \, n \geq 2 .
\end{equation}
Formula (\ref{eqn:4.36}) arises when comparing the various functional
equations satisfied by the series $M_a, R_a$ and $S_a$ -- see Lemma
6.2 of \cite{MN08}.

Now, for every $n \geq 2$ one has a natural identification between 
$NC(n-1)$ and $\{ \alpha \in NC(n) \mid \alpha \ll 1_n \}$. Indeed,
the latter set consists of the partitions $\alpha$ in $NC(n)$ such 
that $1$ and $n$ belong to the same block of $\alpha$; and any such 
partition is uniquely obtained from a $\gamma \in NC(n-1)$ by 
adjoining the number $n$ to the block of $\gamma$ which contains 1.
By using this identification, the summation on the right-hand side 
of (\ref{eqn:4.36}) can be turned into a summation over 
$\{ \alpha \in NC(n) \mid \alpha \ll 1_n \}$. It is suggestive to
write the ensuing equation in the form 
\begin{equation}  \label{eqn:4.37}
\kk_n  = \sum_{ \begin{array}{c}
{\scriptstyle \alpha \in NC(n),}   \\
{\scriptstyle \alpha \ll 1_n}
\end{array} } \ 
\Bigl( \, \prod_{ \begin{array}{c}
{\scriptstyle U \in \alpha,} \\
{\scriptstyle 1_n-special}
\end{array} } \ t_{ |U|-1 } \, \Bigr) \cdot 
\Bigl( \, \prod_{ \begin{array}{c}
{\scriptstyle V \in \alpha,} \\
{\scriptstyle not \ 1_n-special}
\end{array} } \ t_{ |V| } \, \Bigr) ,
\end{equation}
where the product over $U$ has in fact only one factor (a partition
$\alpha \ll 1_n$ has a unique $1_n$-special block, the one which 
contains $1$ and $n$). It is a straightforward exercise, left to 
the reader, to verify that (\ref{eqn:4.37}) can be upgraded to the 
statement that for every $n \geq 1$ and every $\beta \in NC(n)$ 
one has
\begin{equation}  \label{eqn:4.38}
\prod_{W \in \beta} \, \kk_{|W|}  = 
\sum_{ \begin{array}{c}
{\scriptstyle \alpha \in NC(n), }   \\
{\scriptstyle \alpha \ll \beta}
\end{array} } \ 
\Bigl( \, \prod_{ \begin{array}{c}
{\scriptstyle U \in \alpha,} \\
{\scriptstyle \beta-special}
\end{array} } \ t_{ |U|-1 } \, \Bigr) \cdot 
\Bigl( \, \prod_{ \begin{array}{c}
{\scriptstyle V \in \alpha,} \\
{\scriptstyle not \ \beta-special}
\end{array} } \ t_{ |V| } \, \Bigr) .
\end{equation}
Finally, we sum over $\beta \in NC(n)$ on both sides of 
(\ref{eqn:4.38}) and we invoke (\ref{eqn:4.35}) on the left-hand 
side, and (\ref{eqn:4.31}) follows.
\hfill{$\blacksquare$}

\begin{remark}   \label{rem:4.4}
The double sum over $\alpha$ and $\beta$ from Equation 
(\ref{eqn:4.31}) can be treated as an iterated sum in two ways. 
One of them -- sum first over $\alpha$, then over $\beta$ -- has 
in fact just been invoked at the end of the preceding proof. But
actually it is the other order of summation (with $\beta$ first)
which brings the double sum to a simpler form, because it allows
one to take advantage of the Boolean features of the partial 
order $\ll$; specifically, one arrives to Equation (\ref{eqn:4.51})
stated in the next theorem. It is remarkable that (\ref{eqn:4.51}) 
closely resembles the formula used to define the concept of 
``$c$-free cumulants'' in the theory of conditionally free 
convolution -- compare e.g. to the third displayed equation on 
p. 366 of \cite{BLS96}. 
\end{remark}

\begin{theorem}   \label{thm:4.5}
Let $( \cA , \varphi )$ be an algebraic noncommutative probability 
space, let $a \in \cA$ be such that $\varphi (a) = 1$, and consider 
the reciprocal $S$-transform 
$1/S_a (z) = \sum_{n=0}^{\infty} t_n z^n$. Then for every $n \geq 1$
one has
\begin{equation}   \label{eqn:4.51}
\varphi (a^n) = \sum_{ \alpha \in NC(n) } \
\Bigl( \, \prod_{ \begin{array}{c}
{\scriptstyle U \in \alpha,} \\
{\scriptstyle outer}
\end{array} } \ t_{ |U|-1 } \Bigr) \cdot 
\Bigl( \, \prod_{ \begin{array}{c}
{\scriptstyle V \in \alpha,} \\
{\scriptstyle inner}
\end{array} } \ t_{ |V|-1 } + t_{ |V| } \Bigr) .
\end{equation}
\end{theorem}

\begin{proof} 
Let $\alpha$ be a partition in $NC(n)$, and let $\fI , \fO$ 
denote the set of inner blocks and respectively the set of outer 
blocks of $\alpha$. Proposition \ref{prop:2.5} implies that 
\[
\sum_{ \begin{array}{c}
{\scriptstyle \beta \in NC(n),}   \\
{\scriptstyle \beta \gg \alpha}
\end{array} } \ 
\Bigl( \, \prod_{ \begin{array}{c}
{\scriptstyle U \in \alpha,} \\
{\scriptstyle \beta-special}
\end{array} } \ t_{ |U|-1 } \, \Bigr) \cdot 
\Bigl( \, \prod_{ \begin{array}{c}
{\scriptstyle V \in \alpha,} \\
{\scriptstyle not \ \beta-special}
\end{array} } \ t_{ |V| } \, \Bigr) 
\]
\begin{equation}    \label{eqn:4.52}
= \ \sum_{ \begin{array}{c}
{\scriptstyle \fV \subseteq \alpha \ such}   \\
{\scriptstyle that \ \fV \supseteq \fO}
\end{array} } \ 
\Bigl( \, \prod_{U \in \fV} \, t_{ |U|-1 } \, \Bigr) \cdot 
\Bigl( \, \prod_{V \in \alpha \setminus \fV} 
\, t_{ |V| } \, \Bigr) .
\end{equation}
By performing the substitution $\fV \setminus \fO =: \fU$ in the 
sum on the right-hand side of (\ref{eqn:4.52}), this can be 
continued with
\[
= \Bigl( \, \prod_{U \in \fO} \, t_{ |U|-1 } \, \Bigr) \cdot 
\Bigl[ \ \sum_{ \fU \subseteq \fI}
\ \Bigl( \, \prod_{U \in \fU} \, t_{ |U|-1 } \, \Bigr) \cdot 
\Bigl( \, \prod_{V \in \fI \setminus \fU} 
\, t_{ |V| } \, \Bigr) \ \Bigr] .
\]
But it is clear that the latter sum over $\fU \subseteq \fI$ is 
precisely the expansion of the product 
$\prod_{V \in \fI} \bigl( t_{|V|-1} + t_{|V|} \bigr)$. So, 
altogether, what we have obtained is that 
\[
\sum_{ \begin{array}{c}
{\scriptstyle \beta \in NC(n),}   \\
{\scriptstyle \beta \gg \alpha}
\end{array} } \ 
\Bigl( \, \prod_{ \begin{array}{c}
{\scriptstyle U \in \alpha,} \\
{\scriptstyle \beta-special}
\end{array} } \ t_{ |U|-1 } \, \Bigr) \cdot 
\Bigl( \, \prod_{ \begin{array}{c}
{\scriptstyle V \in \alpha,} \\
{\scriptstyle not \ \beta-special}
\end{array} } \ t_{ |V| } \, \Bigr) 
\]
\begin{equation}   \label{eqn:4.53}
= \Bigl( \, \prod_{U \in \fO} \, t_{ |U|-1 } \Bigr) \cdot 
\Bigl( \, \prod_{V \in \fI} \, t_{ |V|-1 } + t_{ |V| } \Bigr) .
\end{equation}

Equation (\ref{eqn:4.53}) holds for every $\alpha \in NC(n)$. Let us
sum over $\alpha$ on both its sides. Then on the left-hand side we 
obtain exactly the double sum over $\alpha$ and $\beta$ which is 
known from Proposition \ref{prop:4.3} to be equal to $\varphi (a^n)$,
and the required formula (\ref{eqn:4.51}) follows.
\end{proof}

$\ $

$\ $

$\ $

Alexandru Nica

Department of Pure Mathematics, University of Waterloo,

Waterloo, Ontario N2L 3G1, Canada.

Email: anica@math.uwaterloo.ca

\end{document}